\documentclass{icmart}

\newtheorem{theorem}{Theorem}[section]

\theoremstyle{definition}

\newcommand{\CC}{\mathbb{C}}

\def\Isom{\mbox{\rm{Isom}}}
\def\span{\mbox{\rm{span}}}

\def\L{\Lambda }

\def\a{\alpha}
\def\b{\beta}

\def\D{\Delta}

\def\G{\Gamma}

\def\l{\lambda}

\def\w{\omega}

\def\PSL{\mbox{\rm{PSL}}}

\def\Isom{\mbox{Isom}}

\def\qed{ $\Box$}

\def\dim{\mbox{\rm{dim}}}

\def\HH{\mathbb{H}}

\def\CC{\mathbb{C}}

\edef\t@mp{\catcode`\noexpand\~=\the\catcode`\~}%
    \catcode`\~=12
    \def\tild@{~}%
    \t@mp

\title{Virtual betti numbers of symmetric spaces} 

\contact[agol@math.uic.edu]{MSCS UIC 322 SEO, m/c 249\\
        851 S. Morgan St.\\
        Chicago, IL 60607-7045}

\begin{document}

\large
\begin{abstract} 
We prove that a closed arithmetic hyperbolic 3-manifold with positive first betti number has
virtually infinite first betti number.  
\end{abstract} 



\maketitle
If $X$ is a manifold, let $\mathcal{H}^{k}(X)$ denote the space of harmonic $k$-forms on $X$. 
If $S$ is a symmetric space, with $G=Isom_{0}(S)$ let $\b_{G}^{k}(S)$ be the dimension
of $G$-invariant harmonic $k$-forms on $S$. Clearly, $\b_{G}^{0}(S)=\b_{G}^{n}(S)=1$.

\begin{theorem} \label{IVB}
Let $S$ be an $n$ dimensional symmetric space, such that  $G = \Isom_{0}(S)$ is simple, and $\L< G$ a uniform arithmetic torsion-free lattice. Suppose
that  $\b^{k}(S/\L) > \b^{k}_{G}(S)$. Then there exists uniform lattices $\L_{i} < \L$ such 
that $\b^{k}(S/\L_{i}) \geq i$.
\end{theorem}
\begin{proof}
Since $\b^{k}(S/\L)>0$, by Hodge's theorem, there exists a non-zero harmonic $k$-form $\w\in \mathcal{H}^{k}(S/\L)$. 
The pullback of $\w$ to $S$ is a harmonic form $\tilde{\w}\in \mathcal{H}^{k}(S)$, such that $\tilde{\w}$ is invariant under
$\L$. That is, for all $g\in \L$,  $g_{\ast}(\tilde{\w})=(g^{-1})^{\ast}(\tilde{\w})=\tilde{\w}$. Since $\b^{k}(S/\L)> \b^{k}_{G}(S)$,
we may assume that $\tilde{\w}$ is not $G$-invariant. 
At each point $x\in S$, there is a canonical norm and inner product on $\L^{k}(T_{x}(S))$, such that the collection
of these norms is equivariant under the action of $G$ on $S$. The space $\mathcal{H}^{k}(S)$ may be given an inner product $\langle,\rangle$, 
by integrating the canonical inner product on $\L^{k}(T_{x}(S))$  over a compact ball $B_{x}(R)$, for
some $x\in S, R>0$. The action of $G$ on $\mathcal{H}^{k}(S)$ is continuous with respect to the topology on $\mathcal{H}^{k}(S)$
induced by this inner product. Since a form $\a\in \mathcal{H}^{k}(S)$ is harmonic, clearly $\|\a\|=0$ if and only if $\a=0$. 
Consider the vector space $V=\span \{ g_{\ast}(\tilde{\w}) | g\in Comm(\L) \}$. Clearly, $V\subset \mathcal{H}^{k}(S)$,
since $\tilde{\w} \in \mathcal{H}^{k}(S)$ and $G$ preserves $\mathcal{H}^{k}(S)$. Also, $V$ is invariant under the action of $Comm(\L)$.  
If $V$ is finite dimensional, then there are elements $g_{1}, \ldots, g_{m} \in Comm(\L)$ such that $V=\span \{ g_{i\ast}(\tilde{\w}), i=1,\ldots, m\}$. 
Let $\G = \cap_{i} (g_{i} \L g_{i}^{-1})$. Since $\tilde{\w}$ is invariant under $\L$, $g_{i\ast}(\tilde{\w})$ is invariant
under the action of $g_{i}\L g_{i}^{-1}$. Thus, $g_{i\ast}(\tilde{\w})$ is invariant under the action of $\G$ for all $i$, and
therefore $\G$ acts trivially on $V$.
We must show that $V\subset \mathcal{H}^{k}(S)$ is a closed subspace. Suppose that 
$v_{\infty}\in \mathcal{H}^{k}(S)$ is a limit of elements $v_{j} \in V$. Since $v_{\infty}\notin V$, 
and $V$ has an orthonormal basis $\{ u_{i}, i=1,\ldots, m\}$, we may take $u=v_{\infty} -\sum_{i=1}^{m} \langle u_{i}, v_{\infty}\rangle u_{i}$.
Then $\langle v_{\infty}, u\rangle \neq 0$, since $u$ is the projection of $v_{\infty}$ onto $V^{\perp}$.
But by continuity of the inner product, $0=\langle v_{j}, u \rangle \to \langle v_{\infty}, u\rangle\neq 0$, a
contradiction. 
Let $g\in G$ be an arbitrary element. Then we want to show that $g_{\ast}(\tilde{\w})\in V$.
Since $\overline{Comm(\L)} = G$, we may find a sequence of elements $g_{i}\in Comm(\L)$, 
such that $\underset{i\to\infty}{\lim} g_{i} = g$. Then clearly $\underset{i\to\infty}{\lim} g_{i\ast}(\tilde{\w}) = g_{\ast}(\tilde{\w})$, both
pointwise and in the topology induced by the inner product $\langle, \rangle$ on $\mathcal{H}^{k}(S)$,
since $G$ acts continuously. 
Since $V$ is closed in $\mathcal{H}^{k}(S)$, we see that $g_{\ast}(\tilde{\w})\in V$. 
But this means that $V$ is
invariant under the action of $G$. Since $\G$ acts trivially on $V$, and $\G$ normally generates
$G$, we see that $G$ acts trivially on $V$. This is a contradiction, since $\tilde{\w}\in V$
is not $G$-invariant.

Thus, we may assume that $V$ is infinite dimensional. For any $i$, choose $g_{1},\ldots, g_{i}\in Comm(\L)$
such that $\{ g_{j\ast}(\tilde{\w}) | j=1,\ldots, i \}$ are linearly independent. Then as before, 
let $\G= \cap_{j} (g_{j}\L g_{j}^{-1})$. Then $\span \{ g_{j\ast}(\tilde{\w}) | j=1,\ldots, i \}$ is fixed
by $\G$, and thus $\b^{k}(S/\G) \geq i$ since $\dim(\mathcal{H}^{k}(S/\G)) \geq i$. 
\end{proof}

\begin{theorem}
Let $M$ be an arithmetic hyperbolic 3-manifold, such that $\b^{1}(M)>0$. Then 
$M$ has infinite virtual positive first betti number. 
\end{theorem}
\begin{proof}
Let $S=\HH^{3}$,  $G=\Isom^{+}(\HH^{3}) =\PSL(2,\CC)$, $\pi_{1}(M)=\L<G$. Then $G$ is simple. 
We also have $\b^{1}_{G}(\HH^{3})=0$. So by theorem \ref{IVB}, for any $i>0$, there exists $\L_{i}<\L$
such that $\b^{1}(\HH^{3}/\L_{i})\geq i$. Thus, $M$ has infinite virtual positive first betti number. 
\end{proof}

{\bf Remark:} The result of Theorem \ref{IVB} is originally due to Borel and Serre in the
case that $\L$ is a congruence arithmetic lattice \cite{Borel76}. 
The argument of Theorem \ref{IVB} may also be used to show that for any $G$-invariant
bundle $E\to S$, and $\Psi$ a $G$-invariant elliptic operator on $E$. Then if there exists a 
$\L$ invariant section $\mu: S\to E$, which is not $G$-invariant, and such that $\Psi\mu=0$,
then there exist $\L_{i}<\L$ such that the dimension of $\L_{i}$-invariant solutions to $\Psi\kappa=0$
is $\geq i$. For example, one can show that if $\l$ is a non-zero eigenvalue of $\D$ on $S/\L$, then 
there is $\L_{i}$ such that the multiplicity of $\l$ as an eigenvalue of $\D$ on $S/\L$ is $\geq i$. Also, Theorem \ref{IVB} ought to 
generalize to non-uniform lattices, but would require an understanding
of the Hodge theory for non-uniform lattices. The extension to
general arithmetic groups follows if one shows that the arithmetic
lattice normally generates $Isom_{0}(S)$.

{\bf Acknowledgements:} We thank Alan Reid, Nicolas Bergeron, and Kevin Whyte for helpful conversations. This paper
was influenced by the appearance of the preprint \cite{CLR06}, and we thank the authors for sharing it
with us. We assumed that the result in this note would be well-known
to experts in arithmetic groups, but the appearance of the 
preprint \cite{Venkataramana06} made it clear that this was false.

\bibliographystyle{../hamsplain.bst}
\bibliography{../refs}

\end{document}